\documentclass[10pt,journal,cspaper,compsoc]{IEEEtran}

\ifCLASSINFOpdf
  \DeclareGraphicsExtensions{.pdf,.jpeg,.png}
\else
  \usepackage[dvips]{graphicx}
  \graphicspath{{../eps/}}
  \DeclareGraphicsExtensions{.eps}
\fi

\usepackage[cmex10]{amsmath}
\usepackage{amssymb,amsfonts}

\usepackage{algorithmicx}
\usepackage[ruled]{algorithm}
\usepackage{algpseudocode}

\usepackage{url}
\usepackage{color}
\usepackage{multirow}

\newcommand{\comment}[1]{} % for invisible comments

\newcommand{\note}[1]{} % for invisible notes
\begin{document}
\title{Some Software Packages for Partial SVD Computation}

\author{Zhouchen~Lin
        % <-this % stops a space
\IEEEcompsocitemizethanks{\IEEEcompsocthanksitem Z. Lin is with
Key Lab. of Machine Perception (MOE), School of EECS, Peking
University, Beijing, China, 100871. E-mail: zlin@pku.edu.cn}}
% <-this % stops a space

\markboth{Manuscript to be enriched}%
{Shell \MakeLowercase{\textit{et al.}}: Bare Demo of IEEEtran.cls for Computer Society Journals}

%\IEEEpubid{\makebox[\columnwidth]{\hfill 0000--0000/00/\$00.00~\copyright~2007 IEEE}%

% PRIOR to the title within the \IEEEcompsoctitleabstractindextext IEEEtran
\IEEEcompsoctitleabstractindextext{%
\begin{abstract}
%\boldmath
This technical report introduces some software packages for partial SVD computation, including optimized PROPACK, modified PROPACK for computing singular values above a threshold and the corresponding singular vectors, and block Lanczos with warm start (BLWS). The current version is preliminary. The details will be enriched soon.
\end{abstract}
% not want either math or citations to appear in the abstract.

% Note that keywords are not normally used for peer review papers.
}

\maketitle
\section{Introduction}
This technical report introduces some software packages for partial SVD computation, including optimized PROPACK, modified PROPACK for computing singular values above a threshold and the corresponding singular vectors, and block Lanczos with warm start (BLWS). The current version is preliminary. The details will be enriched soon.

\section{Optimized PROPACK}
PROPACK \cite{Propack} is nowadays widely used in solving nuclear
norm minimization problems in order to compute the partial SVD. We
optimized PROPACK and obtained 15\% to 20\% speed up. The
optimization is mainly by rewriting the Gram-Schmidt
orthogonalization in the Lanczos procedure. The code is
downloadable at
\url{http://www.cis.pku.edu.cn/faculty/vision/zlin/optPROPACK.zip}

\section{Modified PROPACK}
The current PROPACK \cite{Propack} can only output given
\emph{number} of principal singular values and vectors. However,
when solving nuclear norm minimization problems, we are often
faced with singular value thresholding \cite{Cai2008}, which
requires the principal singular values that are greater than a
given \emph{threshold}. So we modified PROPACK to provide this
functionality. The pseudo code is as in Algorithm~\ref{Algm:PSVDthr} \cite{Chen2010-Thesis}.

\alglanguage{pseudocode}
\begin{algorithm}
\caption{Partial SVD by a Threshold} \label{Algm:PSVDthr}
\begin{algorithmic}[1]
\State Input: Matrix $A$, threshold $svthr$, initial vector $p_0$,
size $K$ of the bidiagonal matrix $B$, $\beta_1=\|p_0\|$,
$u_1=p_0/\beta_1$, $i=1$, $minsv=svthr + 1$,

\While{$minsv \geq svthr$}
\While{$i\leq K$}
     \State $r_i=A^Tu_i-\beta_iv_{i-1}$, $r_i=\mbox{reorth}(r_i)$,
     \State $\alpha_i=\|r_i\|$, $v_i=r_i/\alpha_i$,
     \State $p_i=Av_i-\alpha_i u_i$, $p_i=\mbox{reorth}(p_i)$,
     \State $i\leftarrow i + 1$.
\EndWhile \State Compute the SVD of the bidiagonal matrix $B$ to
obtain the singular values, stored in a vector $s$. \State $neig=$
the number of accurate singular values. \State $minsv=$ the
minimum of those accurate singular values. \State Update $K$ by
\begin{equation*}
K=\left\{
\begin{array}{ll}
\min(K+100,\max(2+K,\\
\mbox{ length}(s > svthr)*K/neig)),&\mbox{if
}neig >0,\\
2*K,&\mbox{otherwise}.
\end{array}
\right.
\end{equation*}

\EndWhile
\end{algorithmic}
\end{algorithm}

The code is downloadable at
\url{http://www.cis.pku.edu.cn/faculty/vision/zlin/ThresholdLANSVD.zip}

\section{Block Lanczos with Warm Start (BLWS)}
When solving nuclear norm minimization problems, we have to solve
the partial SVD multiple times. However, the matrix to compute the
partial SVD only changes slightly over iteration. To utilize this
fact, we propose using the block Lanczos with warm start technique
to cut the computation in each iteration \cite{Lin2010-BLWS}. The
code is downloadable at
\url{http://www.cis.pku.edu.cn/faculty/vision/zlin/BLWS.zip}

The code includes the BLWS technique for eigenvalue decomposition
(BL\_EVD.m) and that for singular value decomposition (BL\_SVD.m),
and also an exemplary usage of BL\_SVD in solving the RPCA problem
\cite{Lin2009-IALM}.

{\small
\bibliographystyle{ieee}
\bibliography{rpca_algorithms}
}

\end{document}